\documentclass[a4paper,english,10pt,twoside,article]{memoir}
\usepackage{babel}
\usepackage{amsmath,amssymb}
\usepackage{mathtools}
\usepackage[thmmarks]{ntheorem}
\usepackage[latin1]{inputenc}
\usepackage[T1]{fontenc}
\usepackage{mathrsfs}
\usepackage{graphicx}
\usepackage{fix-cm}
\usepackage{soul}
\sodef\spred{}{.2em}{.9em plus .4em}{1em plus .1em minus .1em}

\makeevenhead{plain}{\thepage}{}{}
\makeoddhead{plain}{}{}{\thepage}
\makeevenfoot{plain}{}{}{}
\makeoddfoot{plain}{}{}{}

\usepackage[all]{xy}
\SelectTips{cm}{10}

\DeclareMathOperator{\OL}{O\Lambda} 
\DeclareMathOperator{\id}{Id}
\DeclareMathOperator{\M}{M}

\DeclareMathOperator{\aut}{Aut}
\DeclareMathOperator{\mor}{Mor}

\DeclareMathOperator{\Det}{Det}
\DeclareMathOperator{\odet}{Odet}
\DeclareMathOperator{\oDet}{ODet}

\DeclareMathOperator{\sgn}{sgn}

\DeclareMathOperator{\DGer}{Dger}
\DeclareMathOperator{\Dger}{Dger}

\newcommand{\Gl}{\mathrm{Gl}}
\newcommand{\Sl}{\mathrm{Sl}}
\newcommand{\oGl}{\mathrm{OGl}}

\renewcommand{\epsilon}{\varepsilon}
\newcommand{\ua}{{\underline{\text{c}}}}

\newcommand{\cdo}{{{\bullet}}}
\newcommand{\oa}{\overrightarrow{\alpha}}

\newcommand{\LV}{\mc{L}}
\newcommand{\LVo}{\mc{L}_+^*}
\newcommand{\mc}[1]{\mathcal{#1}}
\newcommand{\MnV}{\M_n(\mc{V})}
\newcommand{\MnLV}{\M_n(\LV)}
\newcommand{\MnL}{\Lambda_*}
\newcommand{\GnL}{\Lambda_*}

\newcommand{\oMnL}{\OL_*}
\newcommand{\oGnL}{\OL_*}

\newcommand{\MnB}{\M_n(\mc{B})}

\newcommand{\GnV}{\Gl_n(\mc{V})}
\newcommand{\SlnV}{\Sl_n(\mc{V})}
\newcommand{\SlnLV}{\Sl_n(\LV)}
\newcommand{\GnLV}{\Gl_n(\LV)}
\newcommand{\oGnV}{\oGl_n(\mc{V})}
\newcommand{\oGnLV}{\oGl_n(\LV)}

\newcommand{\N}{\mathbb{N}}
\newcommand{\C}{\mathbb{C}}
\newcommand{\Z}{\mathbb{Z}}
\newcommand{\V}{\mc{V}}
\newcommand{\nn}{\mathbf{n}}
\newcommand{\mm}{\mathbf{m}}

\newcommand{\bcde}{{\beta \gamma\delta\epsilon}}
\newcommand{\acde}{{\alpha\gamma\delta\epsilon}}
\newcommand{\abde}{{\alpha\beta \delta\epsilon}}
\newcommand{\abce}{{\alpha\beta \gamma\epsilon}}
\newcommand{\abcd}{{\alpha\beta \gamma\delta}}

\newcommand{\bcd}{{\beta \gamma\delta}}
\newcommand{\acd}{{\alpha\gamma\delta}}
\newcommand{\abd}{{\alpha\beta \delta}}
\newcommand{\abc}{{\alpha\beta \gamma}}

\newcommand{\ab}{{\alpha\beta}}
\newcommand{\ac}{{\alpha\gamma}}
\newcommand{\ad}{{\alpha\delta}}
\newcommand{\bc}{{\beta\gamma}}
\newcommand{\bd}{{\beta\delta}}
\newcommand{\cd}{{\gamma\delta}}

\newcommand{\de}{{\delta\epsilon}}

\newcommand{\Cec}{$\check{\text{C}}$ec}
\newcommand{\Ce}{\check{C}}

\DeclarePairedDelimiter\absv{\lvert}{\rvert}

\DeclarePairedDelimiter\pare{\lparen}{\rparen}

\theoremsymbol{$\blacklozenge$}
\theorembodyfont{}
\newtheorem{definition}{Definition}[chapter]
\newtheorem{lemma}[definition]{Lemma}
\newtheorem{remark}[definition]{Remark}

\newtheorem{const}[definition]{Construction}

\theoremsymbol{$\square$}
\newtheorem*{bevis}{Proof:}

\begin{document}

\pagestyle{plain}
 
{\LARGE \center Orientations and Connective Structures

on 2-vector Bundles 

$\strut$

\large Thomas Kragh\footnote{The author was funded by the Topology in
  Norway Project, and would like to thank John Rognes and Bjørn Jahren
  for many conversations on the subject.}

}

$\strut$

$\strut$

\noindent \textbf{ABSTRACT}. In \cite{MR2466184} a half magnetic
monopole is discovered and describes an obstruction to
creating a determinant $K(ku) \to ku^*$. In fact it is an obstruction
to creating a determinant gerbe map from $K(ku)$ to $K(\Z,3)$. We
describe this obstruction precisely using monoidal categories and
define the notion of oriented 2-vector bundles, which removes this
obstruction so that we can define a determinant gerbe. We also
generalize Brylinskis notion of a connective structure from
\cite{MR2362847} to 2-vector bundles, in a way compatible with the
determinant gerbe. 

\chapter{Introduction}

In \cite{K-theory/0629} the notion of a charted 2-vector bundle is
defined. This is done such that there is a canonical inclusion of
charted gerbes (essentially a subset of charted 2-vector bundles of
rank 1) into these. They also describe a classifying space
$\absv{B\GnV}$ of equivalence classes of rank $n$ 2-vector bundles,
and this is generalized in \cite{2-cat}. The classifying space of
gerbes is $K(\Z,3)$ (see \cite{MR2362847}), and the inclusion
of gerbes into 2-vector bundles defines a map of classifying spaces
\begin{align} \label{inc}
  K(\Z,3) \to \absv{B\GnV}.
\end{align}
In \cite{MR2466184} it is proven that $\pi_3$ of this map sends the
canonical generator to an element divisible by two (modulo torsion) if
$n$ is large enough. An element which multiplied with 2 mod torsion is
the image of the generator is what they call a half magnetic
monopole. Indeed, this makes sense since a gerbe on $S^3$
representing the canonical generator of $\pi_3(K(\Z,3))$ is a
mathematical model for a magnetic monopole. As in \cite{MR2466184} the
existence of the half magnetic monopole provides an obstruction to
creating a determinant map
\begin{align*}
  \absv{B\GnV} \to B(ku^*) \supset \absv{B\Gl_1(\mc{V})},
\end{align*}
which is the identity on $\absv{B\Gl_1(\mc{V})}$ included into
$\absv{B\GnV}$ in the same way gerbes are
included (block sum with and $n-1$ times $n-1$ identity
matrix). Here $ku^*$ denotes the invertible components of $ku$ with
respect to $\otimes$, i.e. $\{-1,1\}\times BU$. Indeed, such a map
composed with the canonical map 
\begin{align*}
  B(ku^*) \to K(\Z,3)
\end{align*}
would yield a retraction of \eqref{inc}, which is impossible because
the half magnetic monopole in $\pi_3$ should then be sent to an element
which multiplied by 2 is a generator.

In section one we describe this obstruction in the framework of
monoidal categories, and define a natural notion of an orientation on
a 2-vector bundle. We also describe a monoidal category $\oGnV$ such
that $\absv{B\oGnV}$ classifies oriented 2-vector bundles, and we have
a forgetful strict monoidal functor from $\oGnV$ to $\GnV$ inducing
the map of classifying spaces
\begin{align} \label{rect}
  \absv{B\oGnV} \to \absv{B\GnV}.
\end{align}
We then describe the precise obstruction to lifting any map $f \colon
X \to \absv{B\GnV}$ to the oriented ``cover'', as a characteristic
class in $H^3(X,\Z/2\Z)$. Proving that there is a fibration
\begin{align*}
  \absv{B\oGnV} \to \absv{B\GnV} \to K(\Z/2\Z,3).
\end{align*}
We then describe a canonical lift of the inclusion of gerbes and
construct a determinant gerbe functor such that we end up with a
retraction
\begin{align*}
  K(\Z,3) \to \absv{B\oGnV} \to K(\Z,3).
\end{align*}

In \cite{MR2362847} Brylinski defines a connective structures on
gerbes. In section 2 we extend this definition to charted 2-vector
bundles, and prove existence and contractibility of choice. This is
done such that the functors inducing the maps in equation \ref{rect}
takes connective structures to connective structure on charted
bundles.


\chapter{Orientations and Construction of Determinant Gerbe}

Many of the definitions in the following are taken directly from
section 2 and section 3 in \cite{K-theory/0629}. However, some are
taken from \cite{2-cat}, but in the language of monoidal categories
- corresponding to bi-categories with one object.

\begin{definition}
  Let $\Sigma$ be the category with
  \begin{itemize}
  \item{one object $\nn=\{1,\dots,n\}$ for all
      non-negative integers $n\in \N_0$,}
  \item{and morphisms the permutations $\Sigma_n$ of $\nn$.}
  \end{itemize}

  Sum $\oplus$ in $\Sigma$ is defined by disjoint union. More precisely:
  on objects it is standard addition in $\N_0$ and induced morphisms on
  $\nn \mathbf{+} \mm$ is defined by order preservingly identifying
  the first $n$ elements with $\nn$ and the last $m$ elements
  with $\mm$.
  
  Product $\otimes$ in $\Sigma$ is defined by product of sets. More
  precisely: on objects
  it is standard multiplication in $\N_0$ and induced morphisms
  on $\nn\mm$ is defined by identifying the elements in $\nn\mm$
  with the elements in $\nn\times\mm$ using lexicographical
  ordering. I.e. the first $m$ elements in $\nn\mm$ is identified
  with $\{1\}\times \mm$ the next $m$ with $\{2\}\times \mm$ etc.

  These operations are strictly associative and has strict units.
  They are also strictly commutative on the level of objects, but
  not on the induced morphisms. However, choosing the obvious
  permutations as coherency isomorphisms it is well-known that
  we get the structure of a bipermutative category
  (see e.g. \cite{MR0494077}).
\end{definition}

\begin{definition}
  Let $\mc{V}$ be the category with
  \begin{itemize}
  \item{one object $\C^n$ for all
      non-negative integers $n\in \N_0$,}
  \item{and morphisms the linear automorphisms $\Gl_n(\C)$ of 
      $\C^n$.}
  \end{itemize}
  The direct sum functor
  \begin{align*}
    \oplus \colon \mc{V\times V \to V}
  \end{align*}
  is defined by $\C^n\oplus \C^m = \C^{n+m}$ on objects and
  on morphisms by identifying the vector spaces in the standard
  way. The tensor product functor is defined on objects
  by $\C^n \otimes \C^m = \C^{nm}$ and on morphisms by
  using the lexicographical ordering. That is - we identify 
  \begin{align*}
    e_1\otimes e_1' ,\cdots, e_1\otimes e_m', e_2\otimes e_1',\cdots,
    e_n\otimes e_m'
  \end{align*}
  with the standard basis in $\C^{nm}$, where $e_1,\cdots,
  e_n$ and $e_1',\cdots, e_m'$ are the standard bases
  for $\C^n$ and $\C^m$ respectively. As above both operations are
  strictly associative with units. Since the choices involved
  in identifying the bases are the same as the choices made
  for the elements in $\Sigma$, the same permutations viewed as
  matrices may serve as coherency isomorphisms. So again we have
  a bipermutative category, but also a canonical bipermutative functor
  \begin{align*}
    S \colon \Sigma \to \V.
  \end{align*}
\end{definition}

\begin{definition} \label{lv}
  Let $\LV$ be the category with
  \begin{itemize}
  \item{one object $\C_n$ for all
      integers $n\in \Z$,}
  \item{and morphisms the linear automorphisms $\C^*_n = \C^*$.}
  \end{itemize}
  We identify the total space of morphisms with $\Z\times \C^*$, 
  and the direct sum functor is then defined by
  \begin{align*}
    (n,a)\oplus (m,b) = (n+m,ab)
  \end{align*}
  and the tensor functor is defined by
  \begin{align*}
    (n,a)\otimes (m,b) = (nm,a^mb^n).
  \end{align*}
  Both products are strictly associative and commutative. So 
  the coherency isomorphisms could be chosen to be
  identities $(n,1)$. However, for our purpose it turns out that we
  need some of the coherency isomorphisms to be different from the
  identities. More precisely: the coherency twist for the sum 
  \begin{align*}
    \ua \colon \C_n \oplus \C_m \to \C_m \oplus \C_n
  \end{align*}
  is defined to be $(n+m,(-1)^{nm})$ and the twist for the product
  \begin{align*}
    \ua \colon \C_n \otimes \C_m \to \C_m \otimes \C_n
  \end{align*}
  is defined to be $(nm,(-1)^{\frac{n(n-1)m(m-1)}{4}})$. As the following
  lemma shows these choices makes $\LV$ into a bipermutative category.

  It is convenient to introduce the following terminology: a law or rule
  in a monoidal category which regardless of coherency isomorphisms
  holds strictly are called \textbf{weakly strict}. This means that
  the term strict is, as usual, reserved for the laws which have the
  identity as coherency isomorphism, and we see that strict implies
  weakly strict.
\end{definition}

It was noted by John Rognes and it is a curios fact that there are
only two possible $E_\infty$-ring structures on the topological space
$\Z\times BU(1) \simeq \Z\times K(\Z,2)$, and that these arise as the
geometric realization of the category above; but with the two
different choices of coherency isomorphisms: the trivial making all
laws strict and the one we defined.

\begin{lemma}
  The above choice of coherency isomorphisms on the category $\LV$
  makes it bipermutative.
\end{lemma}

\begin{bevis}
  To check that we indeed have a permutative structure on $\LV$ one
  could tediously check all the diagrams in the definition of
  a bipermutative category, but a shorter argument using that we
  know $\Sigma$ to be bipermutative goes as follows.

  Let $\LV_+$ be the full sub-category of $\LV$ defined
  by the non-negatively indexed objects. There is a canonical
  functor $\sgn \colon \Sigma \to \LV_+$ which is the obvious
  bijection on objects and which takes the sign on morphisms. This
  preserves sum and tensor, and sends coherency isomorphisms to the
  signs defined in $\LV$ as coherency. Because it is a bijection on
  objects and the fact that $\Sigma$ is bipermutative makes $\LV_+$
  bipermutative. The coherency sign in $\LV$ for any coherency
  isomorphism only depends on the objects indices modulo $4$. So
  extending to negatively indexed objects by the same formulas will
  still satisfy the necessary equations to be bipermutative.
\end{bevis}

\begin{const}
  Define the functor
  \begin{align*}
    \Lambda \colon \V \to \LV
  \end{align*}
  by $\C^n \mapsto \C_n$ on objects and by taking determinants of
  morphisms. This preserves both sum and product because the
  determinant satisfies
  \begin{align*}
    \det(f\oplus g) = \det(f)\det(g)
  \end{align*}
  and
  \begin{align*}
    \det(f\otimes g) = \det(f)^{\dim(g)}\det(g)^{\dim(f)},
  \end{align*}
  where $\dim(f)$ is the dimension of the underlying vector space.
  The latter can be proved using 
  $(f\otimes g) = (f\otimes \id) \circ (\id \otimes g)$. This explains
  the choice of sum and product in $\LV$ and we may think of $\LV$ as the
  top exterior power of $\V$ extended to negative dimensions. For this
  to be a strict bipermutative functor (or even lax bimonoidal) we
  need that it takes coherency isomorphisms to coherency isomorphisms,
  and this was why we needed the non-trivial signs as coherency isomorphisms
  in definition \ref{lv}. So we have a commutative diagram
  \begin{align*}
    \xymatrix{
      \Sigma \ar[rr]^{S} \ar[dr]^{\sgn} && \V \ar[dl]^{\Lambda} \\
      & \LV 
    }
  \end{align*}
  of bipermutative functors.
\end{const}

We wont use the following explicitly, but it describes very well why
this choice of coherency in $\LV$ is important.

\begin{lemma}
  The induced map on classifying spaces
  \begin{align*}
    \Omega B \absv{\Lambda} \colon ku \to
    \absv{\LV}\simeq \Z \times K(\Z,2)
  \end{align*}
  is the projection to the second Postnikov section in the
  category of $\infty$-loop spaces.
\end{lemma}

\begin{bevis}
  It is an $\infty$-loop map by construction, so all we need to
  check is that it is a $\pi_n$-equivalence for $n\leq 2$.

  The functor $\Lambda$ sends objects $\N_0$ to $\Z$ by the standard
  inclusion. So we need only check that the connectivity of the map is
  at least two on the components corresponding to $n\in \N_0$ for
  large enough $n$. This corresponds to being at least 1-connective on
  the space of automorphisms for large $n$ and the determinant
  \begin{align*}
    \det \colon \Gl_n(\C) \to \C^*
  \end{align*}
  satisfies this.
\end{bevis}

\begin{definition}
  For any bipermutative category $\mc{B}$ define
  $\MnB$ as the category with
  \begin{itemize}
  \item{objects $n$ by $n$ matrices
        $E=(E_{ij})_{i,j=1}^n$ of objects in $\mc{B}$,  and}
  \item{morphisms $n$ by $n$ matrices $\phi = (\phi_{ij})_{i,j=1}^n$
      of morphisms in $\mc{B}$, with the obvious sources and targets.}
  \end{itemize}
  We define a monoidal product on $\MnB$ by
  \begin{align*}
    \cdot \colon \MnB\times\MnB \to \MnB
  \end{align*}
  by standard matrix multiplication formula:
  \begin{align*}
    (E\cdot F)_{ik} = \bigoplus_{j=1}^n (E_{ij}\otimes F_{jk}).
  \end{align*}
  We need not specify parenthesis because $\oplus$ is
  strictly associative. This does, however, not in general
  produce a strictly associative product because this
  would imply both distributive laws in $\mc{B}$ holding strictly.
  But there are obvious coherency isomorphisms induced from the
  coherency isomorphisms in $\mc{B}$ making this a
  monoidal category - with a strict unit.
\end{definition}

We could also define sum of matrices and get a bimonoidal
category, but this is not important in the following. This is because
2-vector bundles will be classified by what we could call the units of
this bimonoidal category and so the product structure is the only
relevant structure.

\begin{const}
  Let $\MnL$ denote the functor induced by $\Lambda$ from
  $\MnV$ to $\MnLV$. This is a strict monoidal functor
  because $\Lambda$ is a bipermutative functor.

  Even though $\LV$ is not equipped with the trivial coherency,
  we may still use the fact that the operations are weakly
  strict to define the symmetric monoidal (with respect to $\oplus$)
  functor
  \begin{align*}
    i \colon \LV \to \LV
  \end{align*}
  by $i(n,a)=(-n,a^{-1})=(-1,1)\otimes(n,a)$. This is a very natural
  choice of inverse to $\oplus$, but be warned: it does not provide
  a coherent choice of inverse in the sense of \cite{MR1030988} when
  passing to the classifying space $\absv{\LV}$. So we cannot
  conclude that the induced $E_\infty$-structure is trivial, which we
  know it is not. We will, however, say that we have a weakly strict
  inverse $i$.

  Let $\det$ be the functor from $\MnLV$ to $\LV$
  given by taking determinant with coefficients in $\LV$.
  Since all commutative, associative and distributive laws
  in $\LV$ are weakly strict and we have an weakly strict
  inverse $i$ to $\oplus$ this is well-defined and sends the
  matrix product to the tensor product in $\LV$. This is
  true independently of the unusual choice of coherency
  isomorphisms in $\LV$ - because it would work with coherencies given
  by identities. However, as the following lemma will show $\det$ is
  not monoidal because of the non-trivial choice of coherency
  isomorphisms.

  Define $\Det = \det \circ \MnL$, again this preserves
  products because both $\det$ and $\MnL$ does so, and
  again the following lemma tells us that it is not monoidal.
\end{const}

As mentioned above 2-vector bundles is related to ``units'' in $\MnV$
and we thus need to define what we mean by this.

\begin{definition}
  Let $\LV^*$ be the full subcategory of invertible objects in $\LV$
  with respect to the product $\otimes$. I.e. using the identification
  in definition \ref{lv} we see
  \begin{align*}
    \mor(\LV^*) = \{\pm 1\} \times \C^*.
  \end{align*}
  This is obviously a permutative category with respect to $\otimes$,
  and since the twist for $\otimes$ on the object pair $-1$ and $-1$
  is not the identity we still retain part of the non-trivial
  coherency structure from $\LV$ in $\LV^*$.

  Also define $\LVo$ to be the sub-category with the single object
  $1$ and $\C^*$ as automorphisms. This is the usual way of
  identifying a group with a category, however, we have also given it
  the canonical permutative structure using that the product is
  Abelian, and it also comes with its inclusion of permutative
  categories into $(\LV,\otimes)$ as the ``positive'' units.
\end{definition}

\begin{definition}
  Let $\GnV$  be the full sub-category of $\MnV$ defined by the
  pre-image of $\LV^*$ using the the functor $\Det$. Also define
  $\GnLV$ by the pre-image of $\LV^*$ using the functor $\det$.

  Similarly we may define $\SlnV$ and $\SlnLV$ using
  $\LVo$ instead of $\LV^*$.
\end{definition}

This definition implies that $\Lambda_*$ maps $\GnV$ to $\GnLV$, and
similarly for the $\Sl_n$'s.

The definition of $\GnV$ is equivalent to the definition in
\cite{K-theory/0629}, because the image object in $\LV$ of
$\Det$ is the determinant of the dimension matrix. The restrictions of
$\det$ and $\Det$ to these sub-categories will also be denoted $\det$
and $\Det$. Objects in $\GnV$ and $\GnLV$ are called weakly invertible
matrices.

\begin{lemma} \label{nonmoi}
  The functors $\Det$ and $\det$ are not monoidal (even on the weakly
  invertible matrices) for $n>1$. More precisely: when evaluated on
  the coherent associativity isomorphisms they produce a sign in
  $\C^*_n$, for some $n$, which in some cases is a minus sign.
\end{lemma}

\begin{remark}
  This is a very important fact and is what turns into the need
  for orientations on 2-vector bundles. It is highly related
  to the Grassmann invariant (see \cite{MR689392}). We plan to describe
  this relation better in \cite{Grassman}.
\end{remark}

\begin{bevis}
  The first statement follows from the second because in $\LV$
  we have strict (not just weakly strict) associativity, and so the
  appearance of a minus sign will imply that the functor does not
  preserve the coherency isomorphisms.

  Since $\Lambda_*$ is injective on objects and strict monoidal
  we only need to find a coherency isomorphism in $\GnLV$
  involving objects in the image of $\Lambda_*$, which is sent
  to minus by $\det$.

  The fact that $\det$ only produces signs in $\C^*$ follows
  because any coherency isomorphism in $\GnLV$ is in each
  entry a coherency isomorphism from $\LV$ which is on the form
  $(n,s)$ with $n\in\Z$ and $s\in \{\pm 1\}$, and taking
  determinant involves $\oplus$ and $\otimes$ which only
  multiplies, divides and takes powers of the last coefficients.

  An easy example of this producing a minus sign for $n=2$ is
  \begin{align*}
    \ua \colon
    \pare*{
      \begin{bmatrix}
        1 & 1 \\
        0  & 1
      \end{bmatrix}\cdot
      \begin{bmatrix}
        0 & 1   \\
        1  & 1     
      \end{bmatrix}
    }\cdot\begin{bmatrix}
      1 & 0   \\
      1  & 1     
    \end{bmatrix} \to
      \begin{bmatrix}
        1 & 1 \\
        0 & 1
      \end{bmatrix}\cdot\pare*{
      \begin{bmatrix}
        0 & 1   \\
        1 & 1     
      \end{bmatrix}
    \cdot\begin{bmatrix}
      1 & 0   \\
      1 & 1     
    \end{bmatrix}},
  \end{align*}
  where $k=\C_k$ in $\LV$. This is the automorphism of the object
  \begin{align*}
    \begin{bmatrix}
      3 & 2   \\
      2 & 1     
    \end{bmatrix}
  \end{align*}
  given by the identity in all but the first entry, where
  it is given by the twist for $1\oplus 1$ on the first two factors plus
  the identity on the last, which is -1.
  \begin{align*}
    \det(\ua) = (1\cdot3-2\cdot2,\frac{1^3(-1)^1}{1^21^2}) = (-1,-1)
    \in \Z\times \C^* = \mor(\LV).
  \end{align*}
  
  This example also works for higher $n$ simply by applying
  block sum with identity matrices. If one would like all determinants
  to take the value $1$ on the objects (i.e. on the first factor
  above) we can take block sum with identity on the two outer matrices
  and take block sum with
  \begin{align*}
    \begin{bmatrix}
      0 & 1   \\
      1 & 0     
    \end{bmatrix}
  \end{align*}
  on the middle factor to produce such an example.
\end{bevis}

We will use this sign to create orientation on 2-vector bundles.

\begin{const} \label{orient}
  The oriented versions of $\GnLV$ ($\GnV$ respectively) are
  denoted $\oGnLV$ ($\oGnV$), and defined to have the
  same objects as the unoriented version but morphisms $(f,s)$
  where $f$ is a morphism in the unoriented version of the
  category and $s\in\{\pm 1\}$.

  Composition and monoidal product is given by composition
  and monoidal product in the original category on the first
  factor and in both cases multiplication on the second factor.

  This would describe a trivial product of the monoidal categories
  $\MnLV$ and the monoidal category with one object and $\Z/2\Z$ as
  automorphisms if we did not ``lift'' the coherency isomorphism (or
  associator) in the following non-trivial way.

  First define the sign $\sgn(\ua)$ of an associator $\ua$ in $\GnLV$
  ($\GnV$) to be  the unique sign such that $\det(\ua)=(\pm 1,
  \sgn(\ua))$. Then the associator in the oriented category is defined
  to be $\ua'=(\ua,\sgn(\ua))$, where $\ua$ is the associator for the
  same objects in the unoriented category.

  These fit into the appropriate commutative diagram, i.e. the
  pentagon relation:
  \begin{align*}
    \xymatrix{
      ((AB)C)D \ar[d]^{\ua'_{AB,C,D}} \ar[rr]^{\ua'_{A,B,C}\cdot \id_{D}} &&
      (A(BC))D \ar[rr]^{\ua'_{A,BC,D}} &&
      A((BC)D) \ar[d]^{\id_A \cdot \ua'_{B,C,D}} \\
      (AB)(CD) \ar[rrrr]^{\ua'_{A,B,CD}} &&&&
      A(B(CD)).
    }
  \end{align*}
  Indeed, this is so because firstly: in the first factor we used the
  coherency from $\GnLV$ ($\GnV$) so on this factor the diagram
  commutes, and secondly: $\det{\id_D}=(\pm 1, 1)$ since $D$ is weakly
  invertible and tensoring with this in $\LV$ preserves any sign in
  the other factor so
  \begin{align} \label{sgnprop}
    \sgn(\ua_{A,B,C}\cdot \id_D) =  \sgn(\ua_{A,B,C}),
  \end{align}
  and since $\det$ is a functor we see that the
  sign around the pentagon must multiply to 1. I.e. they compose to
  the same going from top left to bottom right. The argument is
  identical for $\Det$ replacing $\det$.

  So these are monoidal categories, but moreover we may define functors
  \begin{align*}
    \odet &\colon \oGnLV \to \LV^*, \\  
    \oMnL &\colon \oGnV \to \oGnLV, \\
    \oDet &\colon \oGnV \to \LV^*
  \end{align*}
  defined on morphism by
  \begin{align*}
    \odet(f,s)=\det(f)\otimes (1,s)  \\
    \oMnL(g,s)=(\Lambda_*(g),s) \\
    \oDet(g,s)=\Det(g)\otimes (1,s).
  \end{align*}  
  The result of tensoring with $(1,s)$ is just multiplication with the
  sign $s$ on the morphism in $\LV^*$, which uses that we are in
  $\LV^*$ and not $\LV$. So in fact it is very important that we have
  restricted to the weakly invertible matrices. These oriented versions
  preserve products since the unoriented did and the tensor $\otimes$
  in $\LV$ is weakly strict commutative. The new and useful property
  is that they are in fact strict monoidal because the newly defined
  associators are send to identities (not symmetric monoidal since we
  only fix the associators). Indeed, the signs of the new associators
  are chosen such that they cancel with the sign that made $\det$ and
  $\Det$ not be monoidal.

  There are also canonical strict monoidal functors:
  \begin{align*}
    P_{\GnV} &\colon \oGnV \to \GnV \\
    P_{\GnLV} &\colon \oGnLV \to \GnLV
  \end{align*}
  defined by forgetting the sign.
\end{const}

\begin{remark}
  The composite of the $P$'s with $\det$ and $\Det$ is not the
  oriented functors because if we forget the sign we cannot multiply
  by it. However, the diagram
  \begin{align*}
    \xymatrix{
      \oGnV \ar[r]^{\oGnL} \ar[d]^{P_{\GnV}} & \oGnLV \ar[d]^{P_{\GnLV}} \\
      \GnV \ar[r]^{\GnL} & \GnLV
    }
  \end{align*}
  is obviously a commutative diagram of monoidal functors.
\end{remark}

All of the above categories are smooth in the sense that objects
are discrete, the spaces of morphisms are smooth manifolds, and the
sums and products are on morphisms spaces smooth maps.
We are also in the advantages situation that all the products
we work with are strictly associative and symmetrical on the 
level of objects - meaning that the coherency isomorphisms are
automorphisms. This is so simply because we have only one object in
each isomorphism class.
In \cite{2-cat} the classifying space is defined for any 2-category.
The following is a smooth version of this condensed to our case and
rewritten in the language of monoidal categories. We use the smooth
case only because we later wish to put smooth structures on 2-vector
bundles. In the following we follow the notation for ordered open
coverings in \cite{K-theory/0629}.

\begin{definition} \label{SC2V}
  Let $M$ be a smooth para-compact manifold with smooth partitions of
  unity, let $(\mc{B},\cdot)$ be any smooth monoidal category
  with discrete objects, and let $\mc{(U,J})$ be an ordered open
  cover. A smooth principle $\mc{B}$-bundle is
  \begin{enumerate}[1)]
  \item{
      for each $\alpha<\beta$ in $\mc{J}$ an object
      $E^{\ab}$ in $\mc{B}$, such that for each $\alpha<\beta
      <\gamma$ we have
      \begin{align*}
        E^{\ab} \cdot E^{\bc} = E^{\ac}
      \end{align*}
      on the level of objects, and}
  \item{
      for each $\alpha<\beta<\gamma$ we have smooth maps
      \begin{align*}
        \phi^{\abc} \colon U_{\abc} \to \mor(E^{\ab}\cdot E^{\bc},E^{\ac})
        \qquad(=\aut(E^{\ac})),
      \end{align*}
      called the \textbf{coherency maps} such that}
  \item{the diagram
      \begin{align} \label{cocycle}
        \xymatrix{
          E^{\ab} \cdot \pare*{E^{\bc} \cdot E^{\cd}} 
          \ar[rr]^{\ua^{\abcd}} \ar[d]_{\id \cdot \phi^{\bcd}} &&
          \pare*{E^{\ab} \cdot E^{\bc}}\cdot E^{\cd} 
          \ar[d]^{\phi^{\abc}\cdot \id}\\
          E^{\ab} \cdot E^{\bd} \ar[r]_{\phi^{\abd}} & E^{\ad} &
          E^{\ac} \cdot E^{\cd} \ar[l]^{\phi^{\acd}} 
        }
      \end{align}
      commutes for all points in each quadruple intersection $U_{\abcd}$.}   
  \end{enumerate}
  Here $\ua^{\abcd}$ denotes the associator for the product $\cdot$ in
  $\mc{B}$ related to the two different choices of parenthesis. The
  diagram may be thought of as a cocycle condition.
\end{definition}

\begin{definition}
  Let $n\in \N_0$ be a non-negative integer. A \textbf{smooth charted 
    (oriented) 2-vector bundle} $\mc{E}$ of rank $n$ over $M$ is a
  principal $\GnV$-bundle ($\oGnV$-bundle).  
\end{definition}

This is slightly different than the definition in \cite{K-theory/0629},
but in the unoriented case if we ignore the smoothness condition then
up to the equivalence defined below (taken from \cite{2-cat}) this
provides the same equivalence classes.

\begin{definition}
  Two smooth charted (oriented) 2-vector bundles $\mc{E}_i,i=0,1$ over $X$
  are \textbf{equivalent} if they are cobordant. I.e. there exist a
  smooth charted (oriented) 2-vector bundle $\mc{E}$ over $X\times [0,1]$ such
  that $\mc{E}_{\mid X\times\{t\}}=\mc{E}_t$ for $t=0,1$. Here restriction
  uses restriction of the ordered open cover, which removes
  the sets, and their indices, with empty intersection with $X\times \{i\}$,
  and the smooth coherency maps are assumed to be constant in the
  $t$ direction close to the boundary of $I$ - so as to make composition
  of bordisms well-defined in the smooth category.
\end{definition}

We use this definition, as opposed to the one in \cite{K-theory/0629},
even though it is less explicit, because it is easier to work with.

\begin{definition}
  We say that a smooth charted 2-vector bundle is \textbf{orientable} if it
  is equivalent to a smooth charted 2-vector bundle which is the image
  of a smooth charted oriented 2-vector bundle under the strict monoidal
  functor $P_{\GnV}$.
\end{definition}

Note that the definition of $\mc{B}$-bundle is obviously functorial
with respect to strict monoidal functors, because these preserve
products, associators, and compositions. This even works in the smooth
setting because the functor is in fact also smooth.

\begin{lemma} \label{coco}
  For a smooth charted 2-vector bundle the sign
  \begin{align*}
    \sgn(\ua^{\abcd})    
  \end{align*}
  of the determinants of the associators defines a 3-cocycle in the
  \Cec h complex $\Ce^*(\mc{U},\{\pm 1\})$. The represented class in
  \Cec h cohomology depends only on the equivalence class of the
  2-vector bundle.

  Furthermore, this class is zero if and only if the vector bundle is
  orientable.
\end{lemma}

\begin{bevis}
  This is virtually the same argument as in construction
  \ref{orient}. Again we know from the Lane-Stasheff
  pentagon axiom that
  \begin{align*}
    \ua^{\acde} \circ \ua^{\abce} = \pare*{\ua^{\abcd}\cdot\id_{E^{\de}}} \circ
    \ua^{\abde} \circ \pare*{\id_{E^{\ab}}\cdot\ua^{\bcde}}.
  \end{align*}
  Taking $\Det$, and again using its properties (see equation
  \ref{sgnprop}), we get 
  \begin{align*}
    \sgn\pare*{\ua^{\acde}} \sgn\pare*{\ua^{\abce}} = \sgn\pare*{\ua^{\abcd}} 
    \sgn\pare*{\ua^{\abde}} \sgn\pare*{\ua^{\bcde}},
  \end{align*}
  which is the co-cycle condition. Obviously the associated homology class
  only depends on the equivalence class since the inclusions of
  $X\times \{i\}$ into $X\times [0,1]$ is a homotopy equivalence.

  This class is zero if and only if there is a refinement $\mc{(U,J)}$ of the
  ordered open cover such that we have a chain $\alpha$ in
  $\Ce^2(\mc{U'},\{\pm 1\})$ s.t. $\partial \alpha = \Det(\ua)$, but such 
  a choice exactly corresponds to a lift of the smooth coherency maps
  $\phi^{\abc}$ to $(\phi^{\abc},\alpha)$ in the oriented category,
  such that they satisfy diagram \ref{cocycle} also in the oriented
  category.
\end{bevis}

\begin{remark}
  This last part also tells us how many different choices of
  orientations there are on an orientable 2-vector bundle.
\end{remark}

\begin{definition}
  A charted gerbe is a smooth charted $\LVo$-bundle.
\end{definition}

Note that since $\LVo$ has one object with automorphisms $\C^*$ and is
strictly associative this is the same as having a standard 2-cocycle
with coefficients in $\C^*$, and thus these are classified up to
equivalence by the third cohomology class of the base manifold.

\begin{const}
  We wish to construct an inclusion of gerbes into oriented 2-vector
  bundles, which is a lift of the usual inclusion of gerbes into
  2-vector bundles. We start by describing the usual inclusion on the
  level of categories.
  
  A charted gerbe is the same as a charted $\Sl_1(\LV)$-bundle and
  it is also the same as a charted $\Sl_1(\mc{V})$-bundle. The latter
  seen as just a monoidal category has a natural inclusion into
  $\Gl_1(\mc{V})$, which in turn has a natural inclusion into $\GnV$
  by block sum with the $(n-1)\times (n-1)$ identity matrix on
  objects, and block sum with the identity morphism on this matrix on
  the morphisms. Here we
  call the matrix, which is the unit in the monoidal structure, the
  identity matrix, and of course this has an identity morphism. The
  explicit description of this object is a matrix with the object $\C
  \in \mc{V}$ on the diagonal and the object $\C^0\in \mc{V}$ every
  where else. The latter object has only the identity morphism but
  $\C$ of course has other morphisms, and the identity on the matrix
  is just the identity in each entry.

  The construction of block sum is easily seen to be a strict monoidal
  functor
  \begin{align*}
    \GnV \times \Gl_m(\mc{V}) \to \Gl_{n+m}(\mc{V}).
  \end{align*}
  In fact it is a strict bimonoidal functor on the categories of all
  matrices - not just the weakly invertible, but we wont use this.

  This inclusion of gerbes may be lifted to the oriented category
  $\oGnV$. by using the strict monoidal functor
  \begin{align*}
    \colon \Gl_1(\mc{V}) \to \oGl_1(\mc{V}),
  \end{align*}
  which simply puts the sign $1$ on all morphisms. This is indeed
  strict monoidal since the associator in $\Gl_1(\mc{V})$ is the
  associator for $\otimes$ in $\mc{V}$ on the object $\C$, hence it
  is the identity and has sign 1.

  To generalize the block sum we need to incorporate the sign. So we
  define the functor
  \begin{align*}
    S \colon \oGnV \times \oGl_m(\mc{V}) \to \oGl_{n+m}(\mc{V})
  \end{align*}
  by block sum on objects, block sum on the first factor of the
  morphisms, and by multiplying the signs on the second factor. This is
  obviously a product preserving functor (by ignoring the coherencies
  and using weak strictness on the second factor), and for it to be
  strict monoidal it has to preserve the associators. In the first
  factor of the morphisms this follows because the above unoriented
  block sum is strict monoidal, so we only need to check the sign in
  the second factor. This depends on the fact that we restricted to
  the weakly invertible matrices meaning that 
  \begin{align*}
    \Det(S(\ua_n,\ua_m)) &= \Det(\ua_n)\otimes\Det(\ua_m) =
    (\pm 1,\sgn(\ua_n))\otimes(\pm,\sgn(\ua_m)) = \\ 
    & = (\pm 1,\sgn(\ua_n)\sgn(\ua_m)).
  \end{align*}
  So the sign of the block sum of two associators is the product of
  the signs of the associators, which is precisely what we need. The
  resulting inclusion functor from $\LVo$ to $\oGnV$ will be
  denoted $i$. The unoriented inclusion described above is thus
  the composition of functors $P_{\GnV} \circ i$.
\end{const}

\begin{const} \label{detgerbe}
  In \cite{K-theory/0629} and \cite{2-cat} it is described how to
  construct a simplicial classifying category $B\mc{B}$ of a monoidal
  category $\mc{B}$ such that the geometric realizations of the nerves
  of the categories
  \begin{align*}
    \absv{B\mc{B}}
  \end{align*}
  is a classifying space of $\mc{B}$-bundles. There are certain
  conditions that $\mc{B}$ must satisfy, but the categories we work
  with here satisfy all these. It is used repeatedly in the following 
  that this construction is functorial from the category of monoidal
  categories and strict monoidal functors.

  Since $\LVo$ has one object with automorphisms $\C^*$ it
  follows that $\absv{B\LVo}$ is a $K(\Z,3)$ and in
  \cite{MR2466184} it is proven that the map
  \begin{align*}
    \absv{B(P_{\GnV} \circ i)}_* \colon \pi_3(\absv{B\LVo}) \to
    \pi_3(\absv{B\GnV})
  \end{align*}
  sends the canonical generator to an element divisible by 2 modulo
  torsion. They use this as an obstruction to creating a retraction
  back to $\absv{B\LVo}\simeq K(\Z,3)$. Or more specifically a
  determinant like map to $\absv{B\Gl_1(ku)}=B(\{\-1,1\}\times BU)$,
  which composed with the canonical map to $K(\Z,3)$ would yield such
  a retraction.

  The point of the orientations is that the the monoidal functor $\oDet$
  provides a retraction:
  \begin{align*}
    \absv{B\oDet} \colon \absv{B\oGnV} \to \absv{B\LV^*} \simeq
    K(\Z/2\Z,1)\times \absv{B\LVo} \to \absv{B\LVo}
  \end{align*}
  of $\absv{B i}$. Here the latter map is the projection and the
  identification
  \begin{align*}
    \absv{B\LV^*} \simeq  K(\Z/2\Z,1)\times \absv{B\LVo}
  \end{align*}
  is due to the fact that the monoidal structure in $\LV^*$ is
  strictly associative, and so the $A_\infty$ structure splits. This
  splitting can be described using monoidal functors in the following
  way: we have the inclusion
  \begin{align*}
    \LVo \to \LV^*,
  \end{align*}
  which is a strict symmetric monoidal functor, and we have a
  left inverse (or projection) 
  \begin{align} \label{gerbfunct}
    p \colon \LV^* \to \LVo
  \end{align}
  given by $p(d,a)=(1,a^d)$ this is not symmetric monoidal because the
  symmetry on $-1$ is not the identity. It is, however, strict
  monoidal because
  \begin{align*}
    p((d,a)\otimes(e,b)) &= p(de,a^eb^d) = (1,a^db^e) = \\
    & = (1,a^d)\otimes (1,b^e) = p((d,a))\otimes p((e,b)),
  \end{align*}
  and all associators are identities.

  We now see that $p\circ \oDet$ is a strict monoidal functor and
  right inverse to the strict monoidal functor $i$.
  \begin{definition} \label{detger}
    We define the strict monoidal functor $\DGer$ from $\oGnV$ to
    $\LVo$ as the composite $p \circ \oDet$.
  \end{definition}
\end{const}

We thus conclude that we have removed the obstruction and that the
half magnetic monopole in \cite{MR2466184} must be an unorientable
2-vector bundle.


\chapter{Connective structures}

In \cite{MR2362847} Brylinski defines gerbes over a space $X$ and
classify their equivalence classes by $H^3(X,\Z)$. He also gives an
example in chapter 7 of how to get from a charted gerbe to his
definition of a gerbe, which is easily generalized. He defines the
notion of a connective structure on a gerbe, and in the charted case
this corresponds to having Hermitian connections $\nabla_{\ab}$ on the
trivial line bundles
\begin{align*}
  L_{\ab} = U_{\ab} \times \C
\end{align*}
such that the coherency maps $\phi^{\abc} \colon \C^* \to \C^*$
describes isometries of line bundles
\begin{align*}
  L_{\ab\mid U_\abc} \otimes L_{\bc\mid U_\abc} \to L_{\ac\mid U_\abc}.
\end{align*}

We now generalize this notion of connective structures to 2-vector bundles.

\begin{definition}
  A connective structure $\nabla$ on a smooth charted two vector bundle $\mc{E}$
  is for each $\alpha<\beta$ in $\mc{J}$ a choice of an $n\times n$-matrix
  of connections $\nabla^{\ab}$ on the matrix of trivial bundles
  \begin{align*}
    U_{\ab} \times E^{\ab},
  \end{align*}
  such that for all $\alpha<\beta<\gamma$ the pullback $(\phi^{\abc})^* \nabla^{\ac}$
  is the same connection as the one induced from the product $E^{\ab}\cdot
  E^{\bc}$ of matrices in $\MnV$. This we write as
  \begin{align*}
    (\phi^{\abc})^* \nabla^{\ac} = \nabla^{\ab} \cdot \nabla^{\bc},
  \end{align*}
  and diagram \ref{cocycle} tells us that
  \begin{align*}
    &(\phi^{\acd-1})^* [
    (\phi^{\abc-1})^*(\nabla^{\ab}\cdot \nabla^{\bc})
    \cdot \nabla^{\cd}] = \\
    &(\phi^{\abd-1})^* [\nabla^{\ab} \cdot
    (\phi^{\bcd-1})^*(\nabla^{\bc} \cdot \nabla^{\cd})].
  \end{align*}
  In light of this we define for any such connections  
  the convenient associative ``product''
  \begin{align*}
    \nabla^{\ab}\cdo \nabla^{\bc} = 
    ((\phi^{\abc})^{-1})^* (\nabla^{\ab}\cdot \nabla^{\bc}),
  \end{align*}
  which is a connection on $U_\abc \times E^{\ac}$. In this notation
  the requirement for a family of connections to be a connective
  structure is the cocycle condition with respect to $\cdo$.
\end{definition}

This product behaves well under smooth convex combinations: that is
\begin{align} \label{intercon}
  & (\psi_1 \nabla_1^{\ab} + \psi_2 \nabla_2^{\ab})\cdo
  (\psi'_1 \nabla_1^{\bc} + \psi_2' \nabla_2^{\bc}) = \notag \\
  = & \psi_1\psi_1' \nabla_1^{\ab}\cdo\nabla_1^{\bc} +
  \psi_2\psi_1' \nabla_2^{\ab}\cdo\nabla_1^{\bc} +
  \psi_1\psi_2' \nabla_1^{\ab}\cdo\nabla_2^{\bc} +
  \psi_2\psi_2' \nabla_2^{\ab}\cdo\nabla_2^{\bc},
\end{align}
and
\begin{align} \label{convexcon}
  (\psi_1 \nabla_1^{\ab} + \psi_2 \nabla_2^{\ab})\cdo
  (\psi_1 \nabla_1^{\bc} + \psi_2 \nabla_2^{\bc}) =
  \psi_1 \nabla_1^{\ab}\cdo \nabla_1^{\bc} + 
  \psi_2 \nabla_2^{\ab}\cdo \nabla_2^{\bc}
\end{align}
are true for any smooth functions $\psi_1+\psi_2=1$ and
$\psi_1'+\psi_2'=1$. Indeed, they are true for the tensor product of
vector bundles, and this we may use on each direct summand of each
entry in the matrix, and pullback preserves convex combinations.

\begin{lemma}
  Let $\mc{E}$ be any smooth charted 2-vector bundle. After a possible
  elementary refinement $\mc{E}$ has a connective structure and such
  a choice is a contractible choice.
\end{lemma}

\begin{bevis}
  For any partially ordered set $\mc{J}$ define
  \begin{align*}
    \mc{J}_\alpha^\beta = \{\alpha_0<\cdots<\alpha_k\mid k\in \N, \alpha_0=\alpha,
    \alpha_k=\beta \}
  \end{align*}
  the finite sequences in $\mc{J}$ connecting $\alpha$ and $\beta$.
  By para-compactness we may assume after refinement that
  \begin{itemize}
  \item[C1)]{the cover $(\mc{U},\mc{J})$ is locally finite.}
  \end{itemize}
  Again using para-compactness we shrink $(\mc{U,J})$ to an ordered open
  cover $(\mc{U}',\mc{J})$ (with carrier function the identity),
  such that
  \begin{itemize}
  \item[C2)]{for any $\alpha \in \mc{J}$ we have the
      closure of $U'_\alpha$ contained in $U_\alpha$.}
  \end{itemize}
  Now we use existence of smooth partition of unity to get smooth functions
  \begin{align*}
    \psi_\alpha \colon M \to [0,1]    
  \end{align*}
  with $\psi_{\alpha\mid U'_\alpha}=1$ and $\psi_{\alpha\mid M-U_{\alpha}}=0$.

  Choose any connections $\nabla^{\ab}_0$ on 
  \begin{align*}
    U_\ab \times E^\ab.
  \end{align*}
  For each $\oa\in \mc{J}_\alpha^\beta$ we have a connection
  \begin{align*}
    \nabla^{\oa}_0 = 
    \nabla^{\alpha_0\alpha_1}_{0} \cdo \cdots \cdo 
    \nabla^{\alpha_i\alpha_{i+1}}_{0} \cdo \cdots \cdo 
    \nabla^{\alpha_{k-1}\alpha_k}_{0}
  \end{align*}
  defined on ${U_{\alpha_0\dots\alpha_k}} \subset U_{\ab}$.
  We will smoothly interpolate between these. So we define the
  weights:
  \begin{align*}
    \psi_{\oa} =        
    \smashoperator{\prod_{i=0}^{k-1}} \pare*{ \psi_{\alpha_i}\psi_{\alpha_{i+1}}
      \smashoperator{\prod_{\alpha_i< \gamma <\alpha_{i+1}}}
      (1-\psi_{\gamma})},
  \end{align*}
  for each $\oa\in\mc{J}_\alpha^\beta$. The product is well-defined
  and smooth because C1) implies locally finiteness of the second product.
  Notice that
  \begin{align} \label{prodw}
    \psi_{\oa} = \psi_{\oa_1} \psi_{\oa_2},
  \end{align}
  when $\alpha<\beta<\gamma$, $\oa_1 \in \mc{J}_\alpha^\beta$,
  $\oa_2 \in \mc{J}_\beta^\gamma$, and $\oa \in \mc{J}_\alpha^\gamma$
  is the obvious concatenation of $\oa_1$ and $\oa_2$.  

  On the sets $U'_{\alpha\beta}\subset U_{\alpha\beta}$ we then define
  \begin{align*}
    \nabla^{\ab} = \quad
    \frac{{\sum_{\oa\in\mc{J}_{\alpha}^{\beta}}}
      (\psi_{\oa} \nabla_0^{\oa})}
    {{\sum_{\oa\in\mc{J}_{\alpha}^{\beta}}} (\psi_{\oa})},
  \end{align*}
  This is well-defined and smooth because;
  \begin{itemize}
  \item{the weights in the sum are 0 when the connections are not defined,}
  \item{C1) implies that the sums are locally finite,}
  \item{and $\psi_{\oa}$ is non-zero when we include the $\gamma$'s
      with $\psi_\gamma=1$ and exclude those with $\psi_\gamma=0$ in
      the sequence $\oa$ - notice in particular that $\psi_\alpha=\psi_\beta=1$
      because we restricted to the set $U'_{\ab}$.}
  \end{itemize}

  It satisfies the wanted cocycle condition because: at any point
  $x\in U'_{\abc}$ for $\alpha < \beta < \gamma$ in $\mc{J}$ we have
  $\psi_\alpha=\psi_\beta= \psi_\gamma =1$ and so $\beta$ has to be
  included in the sequence $\oa=(\alpha=\alpha_0 < \cdots <
  \alpha_k=\gamma)$ for $\psi_{\oa}$ to be non-zero. Indeed, we use
  that the subset $\mc{J}_x\subset \mc{J}$ defined by $\mc{J}_x =
  \{\delta \in \mc{J} \mid x\in U_\delta\}$ is totally ordered to
  conclude that if $\beta$ is not in $\oa$ then the factor
  $(1-\psi_\beta)$ is part of the product defining the weight
  $\psi_{\oa}$, which is thus 0. Now we may use \eqref{intercon} and
  \eqref{prodw} repeatedly to see
  \begin{align*}
    \nabla^{\ac}=\frac{
      \smashoperator{\sum_{\mc{J}_\alpha^\gamma}} \psi_{\oa} \nabla_0^{\oa}
    }
    {
      \smashoperator{\sum_{\mc{J}_\alpha^\gamma}} \psi_{\oa}
    }
    = 
    \frac{
      \pare*{
        \smashoperator{\sum_{\mc{J}_\alpha^\beta}} \psi_{\oa} \nabla_0^{\oa}
      }\cdo\pare*{
        \smashoperator{\sum_{\mc{J}_\beta^\gamma}} \psi_{\oa} \nabla_0^{\oa}
      }
    }{
      \pare*{
        \smashoperator{\sum_{\mc{J}_\alpha^\beta}} \psi_{\oa}
      }\pare*{
        \smashoperator{\sum_{\mc{J}_\beta^\gamma}} \psi_{\oa}
      }
    }=\nabla^{\ab}\cdo\nabla^{\bc}.
  \end{align*}
  The choice is contractible because given two connective structures $\nabla$ and
  $\nabla'$ equation \eqref{convexcon} tells us that $t\nabla + (1-t)\nabla'$
  is a connective structure (this convex combination should be interpreted over
  every double intersection).
\end{bevis}

As constructed in the end of the previous section we can use the functor
$\DGer$ to map charted oriented 2-vector bundles to charted gerbes. We
will enrich this map such that it carries connective structures to
connective structures.

\begin{const}
  Let $\mc{E}$ be a charted oriented 2-vector bundle. Taking the
  functor $\Dger$ from definition \ref{detger} on this charted bundle
  produces a 2-cocycle with coefficients in $\C^*$. Indeed, this is
  because the category $\LVo$ to which $\Dger$ maps is the strict
  monoidal category with one object, $\C^*$ as its automorphisms and
  monoidal product the same as composition.

  To get a connective structure we need a connection on each line bundle
  compatible with the isomorphisms. We construct this also using the
  functor $\Dger$: we simply define the parallel transport in 
  $U_{\ab} \times \C$ along a path $f \colon [0,1] \to U_{\ab}$ by
  taking the functor $\Dger$ of the automorphism $P^{\ab}_f$ defined
  on $E^{\ab}$ by parallel transport in
  \begin{align*}
    U_{\ab} \times E^{\ab}
  \end{align*}
  along the path $f$ using the connections $\nabla^{\ab}$. Since
  the connection matrices are preserved using the isomorphisms $\phi^{\abc}$
  we conclude that the parallel transport morphisms are preserved. So
  \begin{align*}
    \xymatrix{
      E^{\ab}\cdot E^{\bc} \ar[rrr]^-{\phi^{\abc}(f(0))}
      \ar[d]_{P^{\ab}_f\cdot P^{\bc}_f}&&&E^{\ac} \ar[d]^{P^{\ac}_f} \\
      E^{\ab}\cdot E^{\bc} \ar[rrr]^-{\phi^{\abc}(f(1))} &&& E^{\ac}
    }
  \end{align*}
  Commutes. This implies that
  \begin{align*}
    \xymatrix{
      \C_1 \ar[rrr]^-{\Dger(\phi^{\abc}(f(0)))}
      \ar[d]_{\Dger(P^{\ab}_f) \otimes \Dger(P^{\bc}_f)}&&& \C_1 \ar[d]^{\Dger(P^{\ac}_f)} \\
      \C_1 \ar[rrr]^-{\Dger(\phi^{\abc}(f(1)))} &&& \C_1
    }
  \end{align*}
  commutes because $\Dger$ is strict monoidal. The monoidal product
  $\otimes$ on $\LVo$ are multiplication in $\C^*$. So the induced
  connections on the line bundles are compatible with the induced
  2-cocycle.
\end{const}  


\bibliographystyle{plain}
\bibliography{/mn/anatu/gjester-u2/thomkr/texmf/tex/Mybib}

\end{document}